\documentclass{article}
\usepackage{amsmath,amscd,amssymb,amsfonts,amsthm}
\usepackage[all]{xy}   
\UseComputerModernTips 
\theoremstyle{plain}
    \newtheorem{theorem}{Theorem}[section]

    \newtheorem{lemma}[theorem]{Lemma}
    \newtheorem{proposition}[theorem]{Proposition}

\theoremstyle{definition}
    \newtheorem{definition}[theorem]{Definition}
    
    \newtheorem{example}[theorem]{Example}
    \newtheorem{remark}[theorem]{Remark}

\DeclareMathOperator{\id}{id}

\def\Z{\mathbb{Z}}

\def\cupMp{\mbox{\raisebox{-1.5ex}{$\stackrel{\sqcup}{\scriptstyle(M,p)}$}}}
\def\cupC{\mbox{\raisebox{-1.5ex}{$\stackrel{\sqcup}{\scriptstyle C}$}}}
\begin{document}
\begin{center}
{\bf \Large Universal localization of triangular matrix rings} \\[1ex]
{\large Desmond Sheiham} {\def\thefootnote{} \footnote{Desmond
Sheiham died on March 25, 2005. This article was prepared for
publication by Andrew Ranicki, with the assistance of Aidan
Schofield.}\addtocounter{footnote}{-1}}
\end{center}
\begin{abstract}
If $R$ is a triangular $2\times 2$ matrix ring, the columns, $P$
and $Q$, are f.g. projective $R$-modules. We describe the
universal localization of $R$ which makes invertible an $R$-module
morphism $\sigma:P\to Q$, generalizing a theorem of A.Schofield.
We also describe the universal localization of $R$-modules.
\end{abstract}
\section{Introduction}
Suppose $R$ is an associative ring (with $1$) and $\sigma:P\to Q$ is a
morphism between finitely generated projective $R$-modules. There is a
universal way to localize $R$ in such a way that $\sigma$ becomes an
isomorphism. More precisely there is a ring
morphism $R\to \sigma^{-1}R$ which is universal for the property that
\begin{equation*}
\sigma^{-1}R\otimes_R P\xrightarrow{1\otimes\sigma} \sigma^{-1}R\otimes_R Q
\end{equation*}
is an isomorphism (Cohn~\cite{Coh71,Coh85,Cohmonthly71,Coh03},
Bergman~\cite{Ber74',Ber78}, Schofield~\cite{Scho85}). Although it is
often difficult to understand universal localizations when $R$ is
non-commutative\footnote{If $R$ is commutative one obtains a ring of
  fractions; see Bergman~\cite[p.68]{Ber78}.} there are examples
where elegant descriptions of $\sigma^{-1}R$ have been possible
(e.g.~Cohn and Dicks~\cite{CohDic76}, Dicks and Sontag~\cite[Thm.
24]{DicSon78}, Farber and Vogel~\cite{FarVog92} ~Ara,
Gonz\'alez-Barroso, Goodearl and Pardo~\cite[Example
2.5]{AGGP04}). The purpose of this note is to describe and to
generalize some particularly interesting examples due to
A.Schofield~\cite[Thm. 13.1]{Scho85} which have application in
topology (e.g.~Ranicki~\cite[Part 2]{Ran03'}).

We consider a triangular matrix ring
$R=\left(\begin{matrix}A & M \\ 0 &
B\end{matrix}\right)$ where $A$ and $B$ are associative rings (with
$1$) and $M$ is an $(A,B)$-bimodule. Multiplication in $R$ is given by
\begin{equation*}
\left(\begin{matrix}
a & m \\ 0 & b
\end{matrix}\right)
\left(\begin{matrix}
a' & m' \\ 0 & b'
\end{matrix}\right)
=
\left(\begin{matrix}
aa' & am'+mb' \\ 0 & bb'
\end{matrix}\right)
\end{equation*}
for all $a,a'\in A$, $m,m'\in M$ and $b,b'\in B$.
The columns
$P=\left(\begin{matrix}A \\
0\end{matrix}\right)$ and $Q=\left(\begin{matrix}M \\
B\end{matrix}\right)$ are f.g. projective left $R$-modules with
\begin{equation*}
P \oplus Q \cong R.
\end{equation*}
The general theory of triangular matrix rings can be found in Haghany
and Varadarajan~\cite{HagVar99,HagVar00}.

We shall describe in Theorem~\ref{structure_of_ring_localization}
the universal localization $R\to\sigma^{-1}R$ which makes
invertible a morphism $\sigma:P \to Q$.
Such a morphism can be written $\sigma=\left(\begin{matrix}j \\
0\end{matrix}\right)$ where $j:A\to M$ is a morphism of left $A$-modules.
Examples follow, in which restrictions are placed on $A$, $B$, $M$ and
$\sigma$. In particular Example~\ref{example:schofields_theorem}
recovers Theorem~13.1 of Schofield~\cite{Scho85}.
We proceed to describe the universal localization
$\sigma^{-1}N=\sigma^{-1}R\otimes_R N$ of an
arbitrary left module $N$ for the triangular matrix ring $R$
(see Theorem~\ref{structure_of_module_localization}).

The structure of this paper is as follows: definitions, statements
of results and examples are given in
Section~\ref{section:results_and_examples} and the proofs are
collected in Section~\ref{section:proofs}.

I am grateful to Andrew Ranicki, Aidan Schofield and Amnon Neeman for
helpful conversations.
\section{Statements and Examples}~\label{section:results_and_examples}
Let us first make more explicit the universal property of localization:
\begin{definition}
A ring morphism $R\to R'$ is called $\sigma$-inverting if
\begin{equation*}
\id\otimes\sigma:R'\otimes_R \left(\begin{matrix}A \\ 0
\end{matrix}\right)\to R'\otimes_R\left(\begin{matrix}M \\ B\end{matrix}\right)
\end{equation*}
is an isomorphism. The universal localization
$i_\sigma:R\to\sigma^{-1}R$ is the initial object in the category
of $\sigma$-inverting ring morphisms $R\to R'$. In other words,
every $\sigma$-inverting ring morphism $R\to R'$ factors uniquely
as a composite $R\to \sigma^{-1}R\to R'$.
\end{definition}
\begin{definition}
An $(A,M,B)$-ring $(S,f_A,f_M,f_B)$ is a ring $S$ together with ring
morphisms $f_A:A\to S$ and $f_B:B\to S$ and an $(A,B)$-bimodule
morphism $f_M:M\to S$.
\begin{equation*}
\xymatrix @=3ex{A\ar[r]^{f_A} & S      & \ar[l]_{f_B} B \\
& M \ar[u]^{f_M} &
}
\end{equation*}
It is understood that the $(A,B)$-bimodule structure on $S$ is induced
by $f_A$ and $f_B$, so that $f_A(a)f_M(m)=f_M(am)$ and
$f_M(m)f_B(b)=f_M(mb)$ for all $a\in A$, $b\in B$ and $m\in M$.

A morphism $(S,f_A,f_M,f_B)\to (S',f_A',f'_M,f'_B)$ of $(A,M,B)$-rings is
a ring morphism $\theta:S\to S'$ such that i) $\theta f_A=f'_A$,
ii) $\theta f_M=f'_M$ and iii) $\theta f_B = f'_B$.
\end{definition}
\begin{definition}
Suppose $p\in M$. Let $(T(M,p),\rho_A,\rho_M,\rho_B)$ denote the
initial object in the subcategory of $(A,M,B)$-rings with the
property $\rho_M(p)=1$. For brevity we often write $T=T(M,p)$.
\end{definition}

The ring $T$ can be explicitly described in terms of generators
and relations as follows. We have one generator $x_m$ for each element
$m\in M$ and relations:
\begin{itemize}
\item[(+)] $x_m+x_{m'} = x_{m+m'}$
\item[(a)] $x_{ap}x_m=x_{am}$
\item[(b)] $x_mx_{pb}=x_{mb}$
\item[($\id$)] $x_p=1$
\end{itemize}
for all $m,m'\in M$, $a\in A$ and $b\in B$. The morphisms $\rho_A$,
$\rho_M$, $\rho_B$ are
\begin{align*}
&\rho_A:A\to T;a\mapsto x_{ap} \\
&\rho_B:B\to T;b\mapsto x_{pb} \\
&\rho_M:M\to T;m\mapsto x_m
\end{align*}

Suppose $\sigma:\left(\begin{matrix}A \\ 0
\end{matrix}\right)\to \left(\begin{matrix}M \\
B\end{matrix}\right)$ is a morphism of left $R$-modules. We may write
$\sigma\left(\begin{matrix}1 \\ 0
\end{matrix}\right) =  \left(\begin{matrix}p \\ 0
\end{matrix}\right)$ for some $p\in M$. Let $T=T(M,p)$.
\begin{theorem}\label{structure_of_ring_localization}
\label{localization_of_matrix_ring}
The universal localization $R\to \sigma^{-1}R$ is (isomorphic to)
\begin{equation*}
R=\left(\begin{matrix}A & M \\ 0 & B\end{matrix}\right)
\xrightarrow{\left(\begin{smallmatrix} \rho_A & \rho_M \\ 0 &
\rho_B\end{smallmatrix}\right)} \left(\begin{matrix} T & T
\\ T & T\end{matrix}\right).
\end{equation*}
\end{theorem}

\begin{example}\label{first_examples}
\begin{enumerate}
\item Suppose $A=B=M$ and multiplication in $A$ defines the $(A,A)$-bimodule
structure on $M$. If $p=1$ then $T=A$ and $\rho_A=\rho_M=\rho_B=\id_A$.
\item Suppose $A=B$ and $M=A\oplus A$ with the obvious bimodule
structure. If $p=(1,0)$ then $T$ is the polynomial ring $A[x]$ in
a central indeterminate $x$. The map $\rho_A=\rho_B$ is the inclusion of
$A$ in $A[x]$ while $\rho_M(1,0)=1$ and $\rho_M(0,1)=x$.
\end{enumerate}
\end{example}
The universal localizations corresponding to Example~\ref{first_examples} are
\begin{enumerate}
\item
$\left(\begin{matrix}
A & A \\
0 & A
\end{matrix}\right)
\to
\left(\begin{matrix}
A & A \\
A & A
\end{matrix}\right)$;
\item
$\left(\begin{matrix}
A & A\oplus A \\
0 & A
\end{matrix}\right)
\to
\left(\begin{matrix}
A[x] & A[x] \\
A[x] & A[x]
\end{matrix}\right)$.
\end{enumerate}
\begin{remark}
One can regard the triangular matrix rings in these examples as path
algebras over $A$ for the quivers
\begin{equation*}
\xymatrix@C=3ex{
1.~~\bullet \ar[r] & \bullet}\qquad 2.~~\xymatrix@=3ex{\bullet \ar@/^/[r] \ar@/_/[r]
& \bullet}
\end{equation*}
The universal localizations $R\to\sigma^{-1}R$ are obtained by
introducing an inverse to the arrow in 1.~and by introducing an
inverse to one of the arrows in 2. See for example
Benson~\cite[p.99]{Ben95} for an introduction to quivers.
\end{remark}
The following examples subsume these:
\begin{example}\label{amal_and_HNN}
\begin{enumerate}
\item (Amalgamated free product; Schofield~\cite[Thm.
4.10]{Scho85}) Suppose $i_A:C\to A$ and $i_B:C\to B$ are ring
morphisms and $M=A\otimes_C B$. If $p=1\otimes 1$ then $T$ is the
amalgamated free product $A \cupC B$ and appears in the pushout
square
\begin{equation*}
\xymatrix@R=3ex{
C \ar[r]^{i_A}\ar[d]_{i_B} & A \ar[d]^{\rho_A} \\
B \ar[r]_(0.4){\rho_B} & T}
\end{equation*}
The map $\rho_M$ is given by $\rho_M(a\otimes b)=\rho_A(a)\rho_B(b)$
for all $a\in A$ and $b\in B$. We recover part 1.~of
Example~\ref{first_examples} by setting $A=B=C$ and $i_A=i_B=\id$.
\item (HNN extension) Suppose $A=B$ and $i_1,i_2:C\to A$ are ring
morphisms. Let $A\otimes_C A$ denote the tensor product with $C$
acting via $i_1$ on the first copy of $A$ and by $i_2$ on the second copy.
Let $M=A\oplus (A\otimes_C A)$ and $p=(1,0\otimes0)$. Now
$T=A*_C\Z[x]$ is generated by the elements in $A$ together with
an indeterminate $x$ and has the relations in $A$ together with
$i_1(c)x=xi_2(c)$ for each $c\in C$. We have $\rho_A(a)=\rho_B(a)=a$
for all $a\in A$ while $\rho_M(1,0\otimes 0)=1$ and $\rho_M(0,a_1\otimes a_2)=
a_1xa_2$. If $C=A$ and $i_1=i_2=\id_A$ we recover part 2.~of
Example~\ref{first_examples}.
\end{enumerate}
\end{example}
The following example is Theorem 13.1 of
Schofield~\cite{Scho85} and generalizes Example~\ref{amal_and_HNN}.
\begin{example}\label{example:schofields_theorem}
\begin{enumerate}
\item
Suppose $p$ generates $M$ as a bimodule, i.e.~$M=ApB$. Now $T$ is
generated by the elements of $A$ and the elements of $B$ subject to
the relation $\sum_{i=1}^n a_ib_i =0$ if $\sum_{i=1}^n a_ipb_i=0$
(with $a_i\in A$ and $b_i\in B$). This ring $T$ is denoted $A\cupMp
B$ in~\cite[Ch13]{Scho85}. The maps $\rho_A$ and $\rho_B$ are obvious and
$\rho_M$ sends $\sum_i a_ipb_i$ to $\sum_i a_ib_i$.
\item Suppose $M$=$ApB\oplus N$ for some $(A,B)$-bimodule $N$. Now
$T$ is the tensor ring over $A\cupMp B$ of
\begin{equation*}
(A\cupMp B)\otimes_A N\otimes_B (A\cupMp B).
\end{equation*}
\end{enumerate}
\end{example}
We may vary the choice of $p$ as the following example illustrates:
\begin{example}\label{first_example_with_interesting_p}
Suppose $A=B=M=\Z$ and $p=2$. In this case $T=\Z\left[\frac{1}{2}\right]$
and $\rho_A=\rho_B$ is the inclusion of $\Z$ in
$\Z\left[\frac{1}{2}\right]$ while $\rho_M(n)=n/2$ for all $n\in\Z$.
\end{example}
Example~\ref{first_example_with_interesting_p} can be verified by
direct calculation using Theorem~\ref{structure_of_ring_localization}
or deduced from part 1.~of Example~\ref{first_examples} by
setting $a_0=b_0=2$ in the following more general proposition. Before
stating it, let us remark that the universal property of $T=T(M,p)$
implies that $T(M,p)$ is functorial in $(M,p)$. An
$(A,B)$-bimodule morphism $\phi:M\to M'$ with $\phi(p)=p'$ induces a
ring morphism $T(M,p)\to T(M',p')$.
\begin{proposition}\label{change_p_to_a0p}
Suppose $A$ and $B$ are rings, $M$ is an $(A,B)$-bimodule and
$p\in M$. If $a_0\in A$ and $b_0\in B$ satisfy $a_0 m=mb_0$ for all
$m\in M$ then
\begin{enumerate}
\item The element $\rho_M(a_0p)=x_{a_0p}=x_{pb_0}$ is central in $T(M,p)$.
\item\label{change_p_localization} The ring morphism $\phi:T(M,p)\to T(M,a_0p)=T(M,pb_0)$ induced by the
bimodule morphism $\phi:M\to M;m\mapsto a_0m=mb_0$ is the universal
localization of $T(M,p)$ making invertible the element
$x_{a_0p}$.
\end{enumerate}
\end{proposition}
\noindent Since $x_{a_0p}$ is central each element in $T(M,a_0p)$
can be written as a fraction $\alpha/\beta$ with numerator
$\alpha\in T(M,p)$ and denominator $\beta=x_{a_0p}^r$ for some
non-negative integer $r$.

Having described universal localization of the ring $R$ in
Theorem~\ref{structure_of_ring_localization} we may also
describe the universal localization $\sigma^{-1}R\otimes_R N$ of a left
$R$-module $N$.
For the convenience of the reader let us first recall the structure of
modules over a triangular matrix ring.
\begin{lemma}\label{describe_R_modules}
Every left $R$-module $N$ can be written canonically as a triple
\begin{equation*}
(N_A,N_B,f:M\otimes_B N_B \to N_A)
\end{equation*}
where $N_A$ is a left $A$-module, $N_B$ is a left $B$-module and
$f$ is a morphism of left $A$-modules.
\end{lemma}
A proof of this lemma is included in Section~\ref{section:proofs}
below. Localization of modules can be expressed as follows:
\begin{theorem}\label{structure_of_module_localization}
    \footnote{This corrects Theorem 2.12 in the preprint version arXiv:math.RA/0407407.} 
For any left $R$-module $N=(N_A,N_B,f)$ the localization left
$\sigma^{-1}R$-module 
$\sigma^{-1}N=\sigma^{-1}R\otimes_R N$ is isomorphic to 
$\left(\begin{matrix} L \\ L \end{matrix}\right)$
with $\sigma^{-1}R=M_2(T)$, $T=T(M,p)$, $L$ the left $T$-module
defined by
$$\begin{array}{ll}
L&=~(T~T)\otimes_RN\\[1ex]
&=~{\rm coker}\big(\begin{pmatrix} 1 \otimes f \\ g \otimes 1
\end{pmatrix}:T\otimes_A M\otimes_B N_B\to
(T\otimes_AN_A) \oplus (T\otimes_B N_B)\big)
\end{array}
$$
with $g$ the $(T,B)$-bimodule morphism
$$g~:~T\otimes_A M \to T~;~t \otimes m \mapsto -t x_m~,$$
and $M_2(T)$ acting on the left of $\begin{pmatrix} L
\\ L \end{pmatrix}$ by matrix multiplication
\end{theorem}

\section{Proofs}\label{section:proofs}
The remainder of this paper is devoted to the proofs of
Theorem~\ref{structure_of_ring_localization},
Proposition~\ref{change_p_to_a0p} and
Theorem~\ref{structure_of_module_localization}.

\subsection{Localization as Pushout}\label{section:Bergman_observation}
Before proving Theorem~~\ref{structure_of_ring_localization} we show
that there is a pushout diagram
\begin{equation*}
\xymatrix{
{\left(\begin{matrix}
\Z & \Z \\
0  & \Z
\end{matrix}\right)} \ar[r]\ar[d]_\alpha &
{\left(\begin{matrix}
\Z & \Z \\
\Z & \Z
\end{matrix}\right)} \ar[d] \\
R \ar[r] & \sigma^{-1}R}
\end{equation*}
where $\alpha\left(\begin{matrix} 1 & 0 \\ 0 & 0
\end{matrix}\right)=\left(\begin{matrix} 1 & 0 \\ 0 & 0
\end{matrix}\right)$, $\alpha\left(\begin{matrix} 0 & 0 \\ 0 & 1
\end{matrix}\right)= \left(\begin{matrix} 0 & 0 \\ 0 & 1
\end{matrix}\right)$ and $\alpha\left(\begin{matrix} 0 & 1 \\ 0 &
0 \end{matrix}\right)=\left(\begin{matrix} 0 & p \\ 0 & 0
\end{matrix}\right)$. Bergman observed~\cite[p.71]{Ber74'}
that more generally, up to Morita equivalence every localization
$R\to \sigma^{-1}R$ appears in such a pushout diagram.

It suffices to check that the lower horizontal arrow in any
pushout
\begin{equation*}
\xymatrix{
{\left(\begin{matrix}
\Z & \Z \\
0  & \Z
\end{matrix}\right)} \ar[r]\ar[d]_\alpha &
{\left(\begin{matrix}
\Z & \Z \\
\Z & \Z
\end{matrix}\right)} \ar[d]_\theta \\
R \ar[r]_i & S}
\end{equation*}
 is i) $\sigma$-inverting and ii) universal among $\sigma$-inverting
 ring morphisms. The universal property of a pushout will be shown to be the universal
 property of a universal localization, so that such a commutative
 diagram is a pushout if and only if $S$ is a universal localization $\sigma^{-1}R$.

i) The map $\id\otimes\sigma : S\otimes_R \left(\begin{matrix}
A \\ 0\end{matrix}\right) \to S\otimes_R \left(\begin{matrix}
M \\ B\end{matrix}\right)$ has inverse given by the composite
\begin{equation*}
S\otimes_R \left(\begin{matrix}
M \\ B\end{matrix}\right) \subset S\otimes_R R \cong S \xrightarrow{\phantom{bbb}\gamma\phantom{bbb}}
S \cong S\otimes_R R \twoheadrightarrow S\otimes_R \left(\begin{matrix}
A \\ 0\end{matrix}\right)
\end{equation*}
where $\gamma$ multiplies on the right by $\theta\left(\begin{matrix}
0 & 0 \\
1 & 0
\end{matrix}\right)$.

ii) If $i':R\to S'$ is a $\sigma$-inverting ring morphism then there
is an inverse $\psi: S'\otimes_R \left(\begin{matrix}
M \\ B\end{matrix}\right) \to S'\otimes_R \left(\begin{matrix}
A \\ 0\end{matrix}\right)$ to $\id\otimes \sigma$. It is argued
shortly below that there is a (unique) diagram
\begin{equation}\label{loc_is_pushout_diagram}
\begin{gathered}
\xymatrix{
{\left(\begin{matrix}
\Z & \Z \\
0  & \Z
\end{matrix}\right)} \ar[r]\ar[d]_\alpha &
{\left(\begin{matrix}
\Z & \Z \\
\Z & \Z
\end{matrix}\right)} \ar[d]_\theta\ar@/^/[ddr]^{\theta'} \\
R \ar[r]^i\ar@/_/[rrd]_{i'} & S\ar@{.>}[dr]  \\
& & S'
}
\end{gathered}
\end{equation}
where $\theta'$ sends $\left(\begin{matrix} 0 & 0 \\ 1 & 0
\end{matrix}\right)$ to $\psi\left( 1 \otimes \left(\begin{matrix}
0 \\ 1\end{matrix}\right)\right)\in S'\otimes_R\left(\begin{matrix}A
\\ 0\end{matrix}\right)\subset S'$. Since $S$ is a pushout there is a
unique morphism $S\to S'$ to complete the diagram and so $i'$ factors
uniquely through~$i$.

To show uniqueness of~(\ref{loc_is_pushout_diagram}),
note that in $S'$ multiplication on the right by
$\theta'\left(\begin{matrix} 0 & 1 \\ 0 & 0
\end{matrix}\right)$ must coincide with the morphism
\begin{equation*}
\left(\begin{matrix} 0 & 0 \\ \id\otimes\sigma & 0
\end{matrix}\right)
\;:\;
S'\otimes \left(\begin{matrix}
A \\ 0\end{matrix}\right) \oplus S'\otimes \left(\begin{matrix}
M \\ B\end{matrix}\right) \longrightarrow S'\otimes \left(\begin{matrix}
A \\ 0\end{matrix}\right) \oplus S'\otimes \left(\begin{matrix}
M \\ B\end{matrix}\right)
\end{equation*}
so multiplication on the right by $\theta'\left(\begin{matrix} 0 & 0 \\ 1 & 0
\end{matrix}\right)$ coincides with $\left(\begin{matrix} 0 & \psi \\ 0 & 0
\end{matrix}\right)$. Now $1\in S'$ may be written
\begin{equation*}
\left(1\otimes\left(\begin{matrix} 1 \\ 0\end{matrix}\right)\ ,\
1\otimes\left(\begin{matrix} 0 \\ 1\end{matrix}\right)\right) \in
S'\otimes_R\left(\begin{matrix} A \\ 0\end{matrix}\right) \oplus S'\otimes_R\left(\begin{matrix} M \\ B\end{matrix}\right)
\end{equation*}
so
$\theta'\left(\begin{matrix} 0 & 0 \\ 1 & 0
\end{matrix}\right) = \psi\left(1\otimes\left(\begin{matrix}0 \\
1\end{matrix}\right)\right)$. The reader may verify that this formula
demonstrates the existence of a commutative diagram~(\ref{loc_is_pushout_diagram}).
\subsection{Identifying $\sigma^{-1}R$}
\begin{proof}[Proof of Theorem~\ref{structure_of_ring_localization}]
It suffices to show that the diagram of ring morphisms
\begin{equation*}
\xymatrix{
{\left(\begin{matrix}
\Z & \Z \\
0  & \Z
\end{matrix}\right)} \ar[r]\ar[d]_\alpha &
{\left(\begin{matrix}
\Z & \Z \\
\Z & \Z
\end{matrix}\right)} \ar[d] \\
{\left(\begin{matrix}
A & M \\
0 & B
\end{matrix}\right)} \ar[r]_\rho & {\left(\begin{matrix}
T & T \\
T & T
\end{matrix}\right)}}
\end{equation*}
is a pushout, where $T=T(M,p)$,  $\rho =
\left(\begin{matrix}\rho_A & \rho_M \\ 0 &
\rho_B\end{matrix}\right)$ and $\alpha$ is defined as in
Section~\ref{section:Bergman_observation}.
Given a diagram of ring morphisms
\begin{equation}\label{pushout_check_diagram}
\begin{gathered}
\xymatrix{
{\left(\begin{matrix}
\Z & \Z \\
0  & \Z
\end{matrix}\right)} \ar[r]\ar[d]_\alpha &
{\left(\begin{matrix}
\Z & \Z \\
\Z & \Z
\end{matrix}\right)} \ar[d]\ar@/^/[ddr]^{\theta} \\
{\left(\begin{matrix}
A & M \\
0 & B
\end{matrix}\right)} \ar@/_/[drr]_{\rho'} \ar[r]_\rho &
{\left(\begin{matrix}
T & T \\
T & T
\end{matrix}\right)}\ar@{.>}[dr]^{\gamma} \\
& & S}
\end{gathered}
\end{equation}
we must show that there is a unique morphism $\gamma$ to complete the diagram.
The map $\theta$ induces a decomposition of $S$ as a matrix ring
$M_2(S')=\left(\begin{matrix} S' & S' \\ S' & S' \end{matrix}\right)$
with $S'$ the centralizer of $\theta(M_2({\mathbb Z})) \subset S$.
In particular, $\theta(e_{ij})=e_{ij}$ for $i,j\in\{1,2\}$.
Any morphism $\gamma$ which makes the
diagram commute must be of the form $\gamma=M_2(\gamma')$ for some
ring morphism $\gamma':T\to S'$ (e.g.~Cohn~\cite[p.1]{Coh85} or
Lam~\cite[(17.7)]{Lam99}). Commutativity of the diagram implies that
$\rho'$ also respects the $2\times2$ matrix structure and we may write
\begin{equation*}
\rho'=\left(\begin{matrix}
\rho'_A & \rho'_M \\
0 & \rho'_B
\end{matrix}\right):
\left(\begin{matrix}
A & M \\
0 & B
\end{matrix}\right)
\longrightarrow
\left(\begin{matrix}
S' & S' \\
S' & S'
\end{matrix}\right)
\end{equation*}
with $\rho'_M(p)=1$ as one sees by considering the images of
$\left(\begin{matrix}0&1\\0&0\end{matrix}\right)$ in
$\left(\begin{matrix}\mathbb{Z}&\mathbb{Z}\\0&\mathbb{Z}
\end{matrix}\right)$
under the maps in the  diagram~(\ref{pushout_check_diagram}) above.
Since $\rho'$ is a ring morphism, one finds
\begin{equation*}
\left(\begin{matrix}
\rho'_A(aa') & \rho'_M(am'+mb') \\
0 & \rho'_B(bb')
\end{matrix}\right)
=
\left(\begin{matrix}
\rho'_A(a)\rho'_A(a') & \rho'_A(a)\rho'_M(m')+\rho'_M(m)\rho'_B(b') \\
0 & \rho'_B(b)\rho'_B(b')
\end{matrix}\right)
\end{equation*}
for all $a,a'\in A$, $b,b'\in B$ and $m,m'\in M$.  Hence the maps
$\rho'_A:A\to S'$ and $\rho'_B:B\to S'$ are ring morphisms and
$\rho'_M$ is a morphism of $(A,B)$-bimodules.  Thus $S'$ is an
$(A,M,B)$-ring with respect to the maps $\rho'_A,\rho'_M,\rho'_B$ such
that $\rho'_M(p)=1$.  By the universal property of $T$ there exists a
unique morphism $\gamma':T\to S'$ such that $M_2(\gamma'):M_2(T)\to
M_2(S')=S$ completes the diagram~(\ref{pushout_check_diagram}) above.
\end{proof}
\begin{proof}[Proof of Proposition~\ref{change_p_to_a0p}]
1. In $T(M,p)$ we have $x_{a_0p}x_m=x_{a_0m}=x_{mb_0}=x_mx_{pb_0}=x_mx_{a_0p}$
   for all $m\in M$. \\ \noindent
2. The map $\phi:M\to M; m\mapsto a_0m$ induces
\begin{align}
\phi:T(M,p)&\to T(M,a_0p) \label{loc_of_T(M,p)} \\
x_m &\mapsto x_{a_0m} \notag
\end{align}
In particular $\phi(x_{a_0p})=x_{a_0^2p}\in T(M,a_0p)$ and we have
\begin{align*}
x_{a_0^2p}x_p = x_{a_0(a_0p)}x_p=x_{a_0p}=1=x_{pb_0} = x_p x_{pb_0^2} =x_px_{a_0^2p}
\end{align*}
so $\phi(x_{a_0p})$ is invertible.

We must check that~(\ref{loc_of_T(M,p)}) is universal. If
$f:T(M,p)\to S$ is a ring morphism and $f(x_{a_0p})$ is
invertible, we claim that there exists a unique $\widetilde{f}:
T(M,a_0p)\to S$ such that $\widetilde{f}\phi=f$. \\ \noindent {\em
Uniqueness}: Suppose $\widetilde{f}\phi=f$. For each $m\in M$ we
have
\begin{equation*}
\widetilde{f}(x_{a_0m})=\widetilde{f}\phi(x_m)=f(x_m).
\end{equation*}
Now
$f(x_{a_0p})\widetilde{f}(x_m)=\widetilde{f}\phi(x_{a_0p})\widetilde{f}(x_m)=
\widetilde{f}(x_{a_0(a_0p)}x_m)=\widetilde{f}(x_{a_0m})=f(x_m)$
so
\begin{equation}\label{define_widetildef}
\widetilde{f}(x_m)=(f(x_{a_0p}))^{-1}f(x_m).
\end{equation}
{\em Existence}: It is straightforward to check that
equation~(\ref{define_widetildef}) provides a definition
of $\widetilde{f}$ which respects the relations
(+),(a),(b) and ($\id$) in $T(M,a_0p)$. Relation (b), for example, is
proved by the equations
\begin{equation*}
\widetilde{f}(x_m)\widetilde{f}(x_{a_0pb})=f(x_{a_0p})^{-1}f(x_m)f(x_{pb})
                                          =f(x_{a_0p})^{-1}f(x_{mb})
                                          =\widetilde{f}(x_{mb})
\end{equation*}
and the other relations are left to the reader.
\end{proof}
\subsection{Module Localization}
We turn finally to the universal localization $\sigma^{-1}R\otimes_R
N$ of an $R$-module $N$.
\begin{proof}[Proof of Lemma~\ref{describe_R_modules}]
If $N$ is a left $R$-module, set $N_A=\left(\begin{matrix} 1 & 0
\\ 0 & 0\end{matrix}\right)N$ and set $N_B=N/N_A$. If $m\in M$ and
$n_B\in N_B$ choose a lift $x\in N$ and define the map $f:M\otimes
N_B\to N_A$ by $f(m\otimes n_B)=\left(\begin{matrix} 0 & m \\ 0 &
0\end{matrix}\right) x$. Conversely, given a triple
$(N_A,N_B,f)$ one recovers a left $R$-module $\left(\begin{matrix} N_A \\
N_B\end{matrix}\right)$ with
\begin{equation*}
\left(\begin{matrix}
a & m \\
0 & b
\end{matrix}\right)
\left(\begin{matrix}
n_A \\
n_B
\end{matrix}\right)
=\left(\begin{matrix}
an_A+f(m\otimes n_B) \\
bn_B
\end{matrix}\right)
\end{equation*}
for all $a\in A$, $b\in B$, $m\in M$, $n_A\in N_A$, $n_B\in N_B$.
\end{proof}
\begin{proof}[Proof of Theorem~\ref{structure_of_module_localization}]
As in the statement, let $T=T(M,p)$ and define the left $T$-module
$$L~=~{\rm coker}\big(\begin{pmatrix} 1 \otimes f \\ g \otimes 1
\end{pmatrix}:T\otimes_A M\otimes_B N_B\to
    (T\otimes_AN_A) \oplus (T\otimes_B N_B)\big)~.$$
We shall establish an isomorphism of left $T$-modules
\begin{equation}\label{simplify_TTtensor}
\left(\begin{matrix} T & T \end{matrix}\right) \otimes_R
\left(\begin{matrix} N_A \\ N_B
\end{matrix}\right) ~\cong~ L
\end{equation}
and leave to the reader the straightforward deduction that there
is an
 isomorphism of $\sigma^{-1}R$-modules
\begin{equation*}
\sigma^{-1}R\otimes_R N =
\left(\begin{matrix}
T & T \\
T & T
\end{matrix}\right)
\otimes_R \left(\begin{matrix}
N_A \\
N_B
\end{matrix}\right)
\cong
\left(\begin{matrix}
L \\
L
\end{matrix}\right).
\end{equation*}
The left $T$-module morphism
$$\begin{array}{l}
\alpha:L \to \left(\begin{matrix} T & T
\end{matrix}\right) \otimes_R \left(\begin{matrix} N_A \\ N_B
\end{matrix}\right)~;\\[1ex]
\hspace*{50pt} (t\otimes n_A,t'\otimes n_B) \mapsto
\left(\begin{matrix}
t & 0 \end{matrix}\right) \otimes_R \left(\begin{matrix} n_A \\
0
\end{matrix}\right) +\left(\begin{matrix} 0 & t' \end{matrix}\right) \otimes_R \left(\begin{matrix} 0 \\
n_B
\end{matrix}\right)
\end{array}$$
is well-defined, since
$$\begin{array}{l}
\alpha(t\otimes_A f(m,n_B),g(t,m)\otimes_B n_B)\\[2ex]
\hspace*{50pt}=~\alpha(t\otimes_A f(m,n_B),-tx_m \otimes_B n_B)\\[2ex]
\hspace*{50pt}=~\begin{pmatrix} t & 0 \end{pmatrix} \otimes_R
\begin{pmatrix} f(m,n_B) \\ 0 \end{pmatrix} -
\begin{pmatrix} 0 & tx_m \end{pmatrix}
\otimes_R \begin{pmatrix} 0 \\ n_B \end{pmatrix}\\[2ex]
\hspace*{50pt}=~\begin{pmatrix} t & 0 \end{pmatrix} \otimes_R
\begin{pmatrix} 0 & m \\ 0 & 0 \end{pmatrix}\begin{pmatrix} 0 \\ n_B \end{pmatrix} -
\begin{pmatrix} t & 0 \end{pmatrix}
\begin{pmatrix} 0 & m \\ 0 & 0 \end{pmatrix}
\otimes_R \begin{pmatrix} 0 \\ n_B \end{pmatrix}\\[2ex]
\hspace*{50pt}=~0 \in \left(\begin{matrix} T & T
\end{matrix}\right) \otimes_R \left(\begin{matrix} N_A \\ N_B
\end{matrix}\right)~.
\end{array}$$
The left $T$-module morphism
$$ \beta: \left(\begin{matrix} T & T \end{matrix}\right) \otimes_R
\left(\begin{matrix} N_A \\ N_B
\end{matrix}\right) \to L~;~\left(\begin{matrix} t & t'
\end{matrix}\right) \otimes_R \left(\begin{matrix} n_A \\ n_B
\end{matrix}\right) \mapsto (t\otimes n_A ,t'\otimes n_B)$$
is well-defined, since
$$\begin{array}{l}
\beta\big(
\begin{pmatrix} t & t'\end{pmatrix}
\otimes_R\begin{pmatrix} a & m
\\ 0 & b \end{pmatrix}
\begin{pmatrix} n_A \\ n_B \end{pmatrix}-
\begin{pmatrix} t & t'\end{pmatrix}
\begin{pmatrix} a & m \\ 0 & b \end{pmatrix} \otimes_R
\begin{pmatrix} n_A \\ n_B \end{pmatrix}\big) \\[2ex]
\hspace*{25pt}=~(t\otimes
(an_A+f(m,n_b)),t'\otimes bn_B)-( ta\otimes n_A,(tx_m+t'b)\otimes n_B)\\[2ex]
\hspace*{25pt}=~(t\otimes f(m,n_B),-tx_m\otimes n_B)\\[2ex]
\hspace*{25pt}=~(1\otimes f,g\otimes 1)(t \otimes m \otimes
n_B)~=~0 \in L~.
\end{array}$$
It is immediate that $\beta\alpha=\id$. To
prove~(\ref{simplify_TTtensor}) we must check that
$\alpha\beta=\id$ or in other words that
\begin{equation*}
\left(\begin{matrix}
t & t'
\end{matrix}\right) \otimes_R
\left(\begin{matrix} n_A \\ n_B
\end{matrix}\right)
=
\left(\begin{matrix}
t & 0
\end{matrix}\right) \otimes_R
\left(\begin{matrix} n_A \\ 0
\end{matrix}\right)
+ \left(\begin{matrix} 0 & t'
\end{matrix}\right) \otimes_R
\left(\begin{matrix} 0 \\ n_B
\end{matrix}\right).
\end{equation*}
This equation follows from the following two calculations:
\begin{align*}
&\left(\begin{matrix}
t & 0
\end{matrix}\right) \otimes_R
\left(\begin{matrix} 0 \\ n_B
\end{matrix}\right) = \left(\begin{matrix}
t & 0
\end{matrix}\right)
\left(\begin{matrix}
0 & 0 \\
0 & 1
\end{matrix}\right) \otimes_R
\left(\begin{matrix} 0 \\ n_B
\end{matrix}\right) =0~; \\
&\left(\begin{matrix}
0 & t'
\end{matrix}\right) \otimes_R
\left(\begin{matrix} n_A \\ 0
\end{matrix}\right) = \left(\begin{matrix}
0 & t'
\end{matrix}\right)
\left(\begin{matrix}
1 & 0 \\
0 & 0
\end{matrix}\right) \otimes_R
\left(\begin{matrix} n_A \\ 0
\end{matrix}\right) =0~.
\qedhere
\end{align*}
\end{proof}

\noindent International University Bremen, Bremen 28759, Germany.
\end{document}